\documentclass[twocolumn,amsmath,amssymb,aps,secnumarabic,%
    nofootinbib,groupedaddress]{revtex4}
\usepackage{times,mathptmx,amsthm,tikz}
\usetikzlibrary{matrix,arrows}

\usepackage[small]{subfigure}

\newtheorem{theorem}{Theorem}[section]
\newtheorem{lemma}[theorem]{Lemma}
\newtheorem{corollary}[theorem]{Corollary}
\theoremstyle{definition}
\theoremstyle{remark}
\newtheorem*{remark}{Remark}
\newtheorem*{example}{Example}

\newcommand{\F}{\mathbb{F}}
\newcommand{\R}{\mathbb{R}}
\newcommand{\C}{\mathbb{C}}
\newcommand{\Z}{\mathbb{Z}}
\newcommand{\Q}{\mathbb{Q}}
\renewcommand{\H}{\mathbb{H}}

\newcommand{\PGL}{\mathrm{PGL}}
\newcommand{\GL}{\mathrm{GL}}
\newcommand{\PSL}{\mathrm{PSL}}
\newcommand{\SL}{\mathrm{SL}}
\newcommand{\U}{\mathrm{U}}
\newcommand{\SU}{\mathrm{SU}}
\newcommand{\PSU}{\mathrm{PSU}}
\newcommand{\PU}{\mathrm{PU}}
\newcommand{\PB}{\mathrm{PB}}
\newcommand{\SO}{\mathrm{SO}}
\renewcommand{\sl}{\mathrm{sl}}
\newcommand{\Tr}{\mathrm{Tr}}

\newcommand{\vc}{\vec{c}}
\newcommand{\cV}{\mathcal{V}}
\newcommand{\mg}{\mathfrak{g}}
\newcommand{\tG}{\tilde{G}}

\newcommand{\eps}{\epsilon}
\renewcommand{\bar}{\overline}

\renewcommand{\tensor}{\otimes}
\newcommand{\braket}[1]{\langle #1\rangle}

\newcommand{\ie}{\textit{i.e.}}

\newcommand{\eq}[2]{\begin{equation}\label{#1}#2\end{equation}}

\newcommand{\thm}[1]{Theorem~\ref{#1}}
\renewcommand{\sec}[1]{Section~\ref{#1}}
\newcommand{\lem}[1]{Lemma~\ref{#1}}
\newcommand{\cor}[1]{Corollary~\ref{#1}}

\newcommand{\fig}[1]{Figure~\ref{#1}}

\newenvironment{tpic}{\begin{tikzpicture}}{\end{tikzpicture}}

\colorlet{darkblue}{blue!70!black}
\colorlet{lightgreen}{green!70!white}
\colorlet{darkred}{red!70!black}

\tikzstyle{tangle}=[draw=white,double=darkblue,line width=.08cm,
    double distance=.6pt]

\tikzstyle{btangle}=[baseline=-.5ex,style=tangle,scale=.25]

\tikzstyle{basic}=[baseline=-.5ex,draw=darkblue,fill=darkblue,semithick,
    scale=.25,>=stealth]

\begin{document}
\title{Denseness and Zariski denseness of Jones braid representations}

\author{Greg Kuperberg}
\email{greg@math.ucdavis.edu}
\thanks{This material is based upon work supported by the National Science
    Foundation under Grant No. 0606795}
\affiliation{Department of Mathematics, University of
    California, Davis, CA 95616}

\begin{abstract}
Using various tools from representation theory and group theory, but without
using hard classification theorems such as the classification of finite
simple groups, we show that the Jones representations of braid groups are
dense in the complex Zariski topology when the parameter $t$ is not a root
of unity. As first established by Freedman, Larsen, and Wang, we the same
result when $t$ is a non-lattice root of unity, other than one initial case
when $t$ has order 10.  We also compute the real Zariski closure of these
representations.  When such a representation is indiscrete in the analytic
topology, then its analytic closure is the same as its real Zariski closure.
\end{abstract}
\maketitle

\section{Introduction}

In this article we will study representations of braid groups associated
with the Jones polynomial $J(L,t)$.  Our question is to determine the
closures of these representations, which are then Lie groups.  Freedman,
Larsen, and Wang \cite{FLW:two} computed these closures in the important
case where $t = \exp(2 \pi i/r)$ is a principal root of unity.  In this
case, the (reduced) braid representations are unitary, and the braids can
be interpreted as quantum circuits.  Freedman, Larsen, and Wang established
that the Jones representations are eventually dense if $r = 5$ or $r \ge 7$.
This has the important corollary that these representations are universal
for quantum computation.

In this article $t$ will usually be a complex number which is not a
root of unity.  Although the Jones polynomial is not directly a model of
quantum computation for these values of $t$, the closure of the braid group
representation is still interesting for related questions in complexity
theory \cite{AAEL:tutte}.  Also, we will say more about the Zariski closure
of the braid group action in the target group $\GL(N,\C)$, rather than
the closure in the usual topology.  Switching to the Zariski topology
simplifies the question, and yet in many cases it does not change the
question very much.  Our main results are as follows:

\begin{theorem} Let $t \in \C$ be a non-zero complex number, and let $n
\ge 4$ and $c \ge 0$ be integers.  Let $X(n \cdot 1,c,t)$ be the (reduced)
Jones representation of the braid group $B_n$ with quantum parameter $t$,
$n$ ordinary strands, and the strand color $c$ at infinity.  If $t$ is not
a root of unity, or if $t$ is a root of unity of order $r \ge 5$ and $r
\ne 6,10$, then the representation is complex Zariski dense in $\SL(X(n
\cdot 1,c,t))$.  When $r = 10$, the same is true with $n+c \ge 5$.
\label{th:complex}
\end{theorem}

A complex number $t$ is a \emph{lattice root of unity} if it is a root of
unity of order 1, 2, 3, 4, or 6.  (So that the ring $\Z[t]$ is a discrete
lattice in $\C$.)  \thm{th:complex} is trivially true when $t = 2$ or $t =
3$, because then the reduced Jones representations are all 1-dimensional.
It is false when $t = 1$ because the Jones representation is large but
trivial; and it is false when $t = 4$ or $t = 6$ for less trivial reasons.
When $r = 10$, the projective image of the braid group action on $X(3 \cdot
1,1,t) = X(4 \cdot 1,0,t)$ lies in $\PSU(2) \cong \SO(3)$ and is that of
the icosahedral group.

\begin{corollary} Assuming the hypotheses of \thm{th:complex}:
\begin{enumerate}
\item If $t \in \C$ and $t \notin \R \cup S^1$ is complex, but neither
real nor norm 1, then the action of $B_n$ is real Zariski dense in $\SL(X(n
\cdot 1,c,t))$.

\item If $t \in \R\setminus \{0,\pm 1\}$, then the action of $B_n$
is real Zariski dense in $\SL(X(n \cdot 1,c,t)_\R)$.

\item If $|t| = 1$ but $t$ is not a root of unity, then $B_n$ acts densely
in $\SU(X(n \cdot 1,c,t))$ in the analytic topology.  Here $\SU(X(n
\cdot 1,c,t))$ is defined using the invariant, but typically indefinite,
Hermitian structure on $X(n \cdot 1,c,t)$.

\item If $t$ is a non-principal non-lattice root of unity, then the action
of $B_n$ is real Zariski dense in $\SU(X(n \cdot 1,c,t))$.

\item If $t$ is a principal non-lattice root of unity, then the action of
$B_n$ is analytically dense in the compact group $\SU(X(n \cdot 1,c,t))$.

\item If the action of $B_n$ on $X(n \cdot 1,c,t)$ is analytically
indiscrete, then in cases 1, 2, and 4 above, the analytic closure is the
same as the real Zariski closure.
\end{enumerate}
\label{c:real} \end{corollary}

\begin{theorem}
Let $c_1, c_2, \ldots, c_\ell$ be distinct non-negative integers and let $t$
and $n$ be as in \thm{th:complex}.  Suppose that for each $k$, the action
of $B_n$ is dense in one of the groups $G_k$ listed in \thm{th:complex} or
\cor{c:real} in either the complex or real Zariski topology. its diagonal
action is dense in the same topology in
$$G_1 \times G_2 \times \cdots \times G_\ell.$$
\label{th:joint} \end{theorem}

Case 5 of \cor{c:real} (plus the analogous case of \thm{th:joint})
is exactly the theorem of Freedman, Larsen, and Wang.  Our argument not
particularly simpler than theirs, but one simplification is that our proof
does not use the classification of finite simple groups.  Also, some of
our ideas and related results were also found by Aharonov, Arad, Eban,
and Landau \cite{AAEL:tutte}.  However, some of our technique, such as
the emphasis on the Zariski topology, is new.

Besides the specific results, our purpose is to describe a set of tools to
establish the closure of representations of groups such as braid groups.
Our tools are described without proof in \sec{s:tools} and established in
\sec{s:toproofs}.  We think that they could be used to further investigate
the images of braid group actions, such as those coming from the colored
Jones polynomial or other quantum link invariants.  (The Jones polynomial
is associated to defining representation of the Lie algebra $\sl(2)$.
There is a polynomial invariant for every simple Lie algebra, colored by
its irreducible representations.)

As an example use of our tools, case 4 of \cor{c:real} follows from case
5, the FLW theorem, using the fact that the complex Zariski topology is
preserved by Galois automorphisms of $\C$.  The action in case 4 is then
dense if and only if there is a braid whose action is an elliptic element
of infinite order.

In \sec{s:other}, we will discuss examples where the braid group acts
discretely or indiscretely on $X(3\cdot 1,1,t) = X(4\cdot 1,0,t)$ and
$X(4\cdot 1,2,t)$ in the analytic topology.  Again, case 6 of \cor{c:real}
says that if $B_n$ acts indiscretely, then it is dense in its real Zariski
closure.

\begin{theorem} The following is a complete classification of when $B_3$
acts discretely on $X(3\cdot 1,1,t)$, or equivalently when $B_4$ acts
discretely on $X(4 \cdot 1,0,t)$, in the cases $t \in \R$ and $|t| = 1$.
\begin{enumerate}
\item Let $t \ne 0,-1$ be real.  Then
$B_3$ acts discretely in $\PGL(X(3\cdot 1,1,t)_\R)$ in the following cases:
$$t < 0 \qquad t + t^{-1} \ge 3 \qquad t + t^{-1} = 1 + 2(\cos \frac{2\pi}{n}).$$

\item Let $t = \exp(i\theta)$ with $0 < |\theta| < 2\pi/3$.  Then $B_3$ acts
discretely in the compact group $\PSU(X(3 \cdot 1,1,t))$ if and only if
$|\theta| = \pi - (2\pi/n)$.

\item Let $t = \exp(i\theta)$ with $\pi > |\theta| > 2\pi/3$.
Then $B_3$ acts discretely in the non-compact group $\PSU(X(3 \cdot 1,1,t))$
if and only if $|\theta| = \pi - (2\pi/n)$.

\item If
$$t + t^{-1} = 1 + 2\cos \frac{2\pi}{7},$$
then $B_4$ acts indiscretely in $\PGL(X(4,2,t)_\R)$.
\end{enumerate}
\label{th:discrete} \end{theorem}

Finally in \sec{s:other}, we will discuss the following corollary of
\thm{th:complex}.  It is used in a previous paper \cite{Kuperberg:jones} to
establish that certain values of the Tutte polynomial are $\mathsf{\#P}$-hard
to estimate multiplicatively.

\begin{corollary} The edge operators $A_{j,y}$ and $B_{j,x}$ defined in
\cite{Kuperberg:jones}, using all $x,y \ne 1$, generate
$\PSL(V(n)_\R)$, where $V(n)_\R$ is the real skein space of the Tutte
polynomial on $n$ vertices with Potts model parameter $q \ge 4$.
\label{c:tutte} \end{corollary}

\acknowledgments

The author would like to thank Scott Carnahan, Nathan Dunfield, Mike
Freedman, Jim Humphreys, Misha Kapovich, Michael Larsen, and David Speyer
for useful discussions.

\section{Definitions for the main results}
\label{s:defs}

\subsection{Geometry}

We the definitions needed in the statements of the main results:  If $\F$
is an infinite field, the \emph{Zariski topology} on the vector space $\F^N$,
or any subset of $\F^N$, is the coarsest topology in which the solution set
of a polynomial equation is closed.  If we interpret $\SL(N,\C)$ as a subset
of $\C^{N^2}$, it inherits a \emph{complex Zariski topology}.  If instead
we interpret it as a subset of $\R^{2N^2}$, it inherits a \emph{real Zariski
topology}. If we further restrict to $\SU(P,Q) \subset \SL(N,\C)$ with $P+Q =
N$, then the real and complex Zariski topologies agree.  Likewise $\SL(N,\R)$
has the same Zariski topology, whether it is viewed as a subset of $\R^{N^2}$
or as a subset of $\SL(N,\C)$ with either of its Zariski topologies.

The analytic topology on any subset of $\R^N$ or $\C^N$ is the usual
topology used in calculus and in most mathematics.

If $G$ is a topological group with a closed subgroup $H$ and another subgroup
$\Gamma$, then we say that $\Gamma$ is dense in $H$ if $\Gamma \cap H$ is
dense in $H$.

An operator $x$ on a vector space $X$ is \emph{elliptic} if it is
diagonalizable, and if its eigenvalues are all on the unit circle.

\subsection{Quantum algebra}

The Jones polynomial can be defined from the Kauffman bracket,
which is defined by these skein relations:
$$\begin{tpic}[style=btangle]
\draw[double] (-1,1) -- (1,-1); \draw[double] (-1,-1) -- (1,1);
\end{tpic}
=
-t^{1/4}
\begin{tpic}[style=basic]
\draw (-1,-1) arc (-45:45:1.414); \draw (1,1) arc (135:225:1.414);
\end{tpic}
-t^{-1/4}
\begin{tpic}[style=basic]
\draw (1,-1) arc (45:135:1.414); \draw (-1,1) arc (225:315:1.414);
\end{tpic}
\quad \qquad
\begin{tpic}[style=basic]
\draw (0,0) circle (1);
\end{tpic}
= -t^{1/2}-t^{-1/2}.$$
For the moment, we let $t^{1/4} \in \F \setminus\{0\}$ for a field $\F$
of characteristic 0.  Note also that the Kauffman bracket is only invariant
under two of the three Reidemeister moves; it gains a factor of $t^{3/4}$
under the first Reidemeister move.

Given non-negative integers $n$ and $c$ of the same parity, define the
\emph{skein space} $W(n \cdot 1,c,t)$ to be the vector space of formal linear
combinations of planar matchings in a rectangle, with $n$ points on the left
and $c$ points on the right.  The elements of $W(n \cdot 1,c,t)$ are \emph{skeins}.
The $c$ points on the right are together called a \emph{clasp}.  We set
to $0$ those matchings that have a U-turn at the right side:
$$\begin{tpic}[style=basic,scale=1.5]
\draw (-3,-2) .. controls (-1,-2) and (1,-1) .. (3,-1);
\draw (-3,-1) .. controls (-1,-1) and (-1,0) .. (-3,0);
\draw (-3,1) .. controls (-1,1) and (1,0) .. (3,0);
\draw (-3,2) .. controls (-1,2) and (1,1) .. (3,1);
\draw[darkgray,dashed] (-3,-3) rectangle (3,3);
\draw[thick,darkred] (3,-1.5) -- (3,1.5);
\end{tpic} \ne 0 \qquad \qquad
\begin{tpic}[style=basic,scale=1.5]
\draw (-3,-2) .. controls (-1,-2) and (1,-1) .. (3,-1);
\draw (-3,-1) .. controls (1,-1) and (1,2) .. (-3,2);
\draw (-3,0) .. controls (-1,0) and (-1,1) .. (-3,1);
\draw (3,0) .. controls (1,0) and (1,1) .. (3,1);
\draw[darkgray,dashed] (-3,-3) rectangle (3,3);
\draw[thick,darkred] (3,-1.5) -- (3,1.5);
\end{tpic} = 0.$$

The braid group $B_n$ acts on $W(n \cdot 1,c,t)$ by braiding the clasps
and then expanding crossings:
$$\begin{tpic}[style=basic,scale=1.5]
\draw (-5,-2) -- (-3,-2) .. controls (-1,-2) and (1,-1) .. (3,-1);
\draw (-5,-1) -- (-3,-1) .. controls (-1,-1) and (-1,0) .. (-3,0) -- (-5,0);
\draw[style=tangle] (-5,2) .. controls (-4,2) and (-4,1) .. (-3,1)
    .. controls (-1,1) and (1,0) .. (3,0);
\draw[style=tangle] (-5,1) .. controls (-4,1) and (-4,2) .. (-3,2)
    .. controls (-1,2) and (1,1) .. (3,1);
\draw[darkgray,dashed] (-3,-3) rectangle (3,3);
\draw[thick,darkred] (3,-1.5) -- (3,1.5);
\end{tpic} = - t^{-1/4}\;\;
\begin{tpic}[style=basic,scale=1.5]
\draw (-3,-2) .. controls (-1,-2) and (1,-1) .. (3,-1);
\draw (-3,-1) .. controls (-1,-1) and (-1,0) .. (-3,0);
\draw (-3,1) .. controls (-1,1) and (1,0) .. (3,0);
\draw (-3,2) .. controls (-1,2) and (1,1) .. (3,1);
\draw[darkgray,dashed] (-3,-3) rectangle (3,3);
\draw[thick,darkred] (3,-1.5) -- (3,1.5);
\end{tpic}\;.$$
When $t$ is not a root of unity, we let $X(n \cdot 1,c,t) = W(n \cdot
1,c,t)$; this skein space with the action of $B_n$ is then projectively
equivalent to the Jones representation.  Projective equivalence does not
change the dense properties of interest to us; see \lem{l:comm}.

We would like to define $W(n \cdot 1,c,t)_\F$ over a field $\F$ that
does not contain $t^{1/4}$ or even $t^{1/2}$.  To this end,
we first let each braid generator be a right-handed crossing times $t^{1/4}$:
$$\tau_2 = t^{1/4}\;\begin{tpic}[style=basic]
\draw[style=tangle] (-1,-1) -- (1,1);
\draw[style=tangle] (-1,1) -- (1,-1); 
\draw (-1,2) -- (1,2); \draw (-1,-2) -- (1,-2);
\end{tpic}\;.$$
Second, we bicolor the complementary regions of a matching
in the rectangle, so that regions are alternately black
and white and the top region is white:
$$
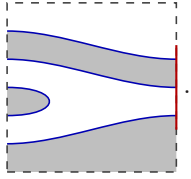
\begin{tpic}[style=basic,scale=1.5]
\begin{scope}[lightgray]
\fill (-3,-3) -- (-3,-2) .. controls (-1,-2) and (1,-1) .. (3,-1)
    -- (3,-3) -- cycle;
\fill (-3,-1) .. controls (-1,-1) and (-1,0) .. (-3,0) -- cycle;
\fill (-3,2) .. controls (-1,2) and (1,1) .. (3,1) --
    (3,0) .. controls (1,0) and (-1,1) .. (-3,1) -- cycle;
\end{scope}
\draw (-3,-2) .. controls (-1,-2) and (1,-1) .. (3,-1);
\draw (-3,-1) .. controls (-1,-1) and (-1,0) .. (-3,0);
\draw (-3,1) .. controls (-1,1) and (1,0) .. (3,0);
\draw (-3,2) .. controls (-1,2) and (1,1) .. (3,1);
\draw[darkgray,dashed] (-3,-3) rectangle (3,3);
\draw[thick,darkred] (3,-1.5) -- (3,1.5);
\end{tpic}\;.$$
Then we can define a matching to be even or odd according to the number
of black regions in the coloring.  If we multiply the odd matchings by
$t^{1/2}$, the result is a basis for which action of $B_n$ is
defined over $\F$.

The skein space $W(n \cdot 1,c,t)$ carries a natural bilinear form which
is used to define $X(n \cdot 1,c,t)$ when $t$ is a root of unity.  If $r$
is the order of $t$, we assume that $c \le r-2$.  Given two skeins $A$
and $B$, we pair them by sewing their rectangles together at both ends.
Where the clasps meet, we insert a special skein called a \emph{Jones-Wenzl
projector} (defined in \sec{s:quantum}):
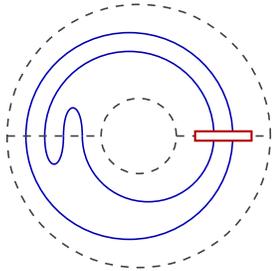
$$\begin{tpic}[style=basic]
\draw (4,0) arc (0:180:4.5);
\draw (5,0) arc (0:180:5.5);
\draw (-3,0) .. controls (-3,2) and (-4,2) .. (-4,0);
\draw (-4,0) .. controls (-4,-2) and (-5,-2) .. (-5,0);
\draw (4,0) arc (0:-180:3.5);
\draw (5,0) arc (0:-180:5.5);
\begin{scope}[dashed,darkgray]
\draw (0,0) circle (2); \draw (0,0) circle (7);
\draw (-7,0) -- (-2,0); \draw (2,0) -- (7,0);
\end{scope}\draw[thick,darkred,fill=white] (3,-.25) rectangle (6,.25);
\end{tpic}\;.$$
The Kauffman bracket value of this diagram is then the value of the pairing.
When $t$ is a root of unity, the form $\braket{A,B}$ has a kernel, and
$X(n \cdot 1,c,t)$ is defined by annihilating this kernel so that $\langle A,B
\rangle$ is non-degenerate in the quotient.

If a skein $A$ has crossings, then we define its reflection $\bar{A}$
by switching left and right crossings.  The form $\braket{\bar{A},B}$
is invariant under the action of $B_n$, and when $|t|=1$ it is Hermitian.

\section{The tools}
\label{s:tools}

\begin{lemma}[Adjoint] Let $\Gamma$ be a subgroup of a simple Lie group $G$.
Then $\Gamma$ is dense if and only if the Lie algebra $\mg$ is irreducible
under the adjoint action of $\Gamma$ and $\Gamma$ is indiscrete.
\label{l:adj} \end{lemma}

\lem{l:adj} reduces the question of computing the closure of a group in a
Lie group to representation theory.  The lemma holds in both the analytic
topology on $G$ and in the Zariski topology, defined below.  Henceforth
we will assume that all representations and other vector spaces that 
we use are finite-dimensional.

In light of \lem{l:adj}, we define a subgroup $\Gamma \subset G$
of a simple Lie group $G$ to be \emph{adjoint-irreducible} if $\mg$
is $\Gamma$-irreducible, but $\Gamma$ is not necessarily dense.  Thus,
\lem{l:adj} by itself establishes case 6 of \cor{c:real}:  If $\Gamma$ (or an
action of $\Gamma$) is real Zariski dense, then it is adjoint-irreducible.
If in addition it is indiscrete in the analytic topology, it is then dense
in that topology.

\begin{lemma}[Zariski] If $\Gamma$ is a Zariski-dense subgroup of a
group $G$, then every admissible representation $V$ of $G$ decomposes in
the same way as a $G$-representation and as a $\Gamma$-representation.
\label{l:zariski} \end{lemma}

\lem{l:zariski} uses the Zariski topology as a convenient converse,
using closures of groups to help compute decompositions of representations.
In order to define the Zariski topology on a group $G$, we first choose
a class $\cV$ of admissible representations over an infinite field $\F$.
The $\cV$ should be closed under passing to duals, direct sums, tensor
products, and subrepresentations (because taking this closure does not
change the resulting topology).  Given $\cV$, the \emph{Zariski topology}
on $G$ is the topology whose closed sets are generated by solutions to
polynomial equations in the entries of the matrices of representations $V
\in \cV$, using some basis for $V$.

If $G \subseteq \GL(N,\F)$, then one standard admissible class is the set
of representations generated from the defining representation $V = \F^N$.
In particular, this $V$ generates the adjoint representation.  We are also
interested in the real Zariski topology on $\GL(N,\C)$.  This topology is
generated by the same representation $V = \C^N$, but interpreted as a real
vector space $V \cong \R^{2N}$ instead.

\begin{lemma}[Compact] A subgroup $G \subseteq \U(N)$ is analytically
closed if and only if it is Zariski closed.
\label{l:compact} \end{lemma}

\lem{l:compact} establishes case 5 of \cor{c:real} from \thm{th:complex}.

\begin{lemma}[Closure] If $f:G \to H$ is an algebraic homomorphism between
algebraic groups, then the image $f(G)$ is Zariski closed.
\label{l:closure} \end{lemma}

An \emph{algebraic group} is an affine algebraic variety that is also a
group with an algebraic group law.  Not every group with a Zariski topology
in our sense is an algebraic group.  However, any Zariski-closed subgroup
of $\GL(N,\F)$ is an algebraic group, and this covers all cases of interest
to us.

\lem{l:closure} is a major reason to trust the Zariski topology in
our context.  Indeed, the Zariski topology is more trustworthy than the
analytic topology, in the sense that \lem{l:closure} is not always true
in the analytic topology for non-compact groups.

\begin{example} Consider the ring $\Z[\sqrt{3}]$ (or the ring of integers
of any other real Galois number field larger than $\Q$).  Let $\Gamma$
be the group of pairs of matrices $(x,\bar{x})$ in
$$G = \SL(2,\R) \times \SL(2,\R),$$
where $x \in \SL(2,\Z[\sqrt{3}])$ and $\bar{x}$ is its Galois conjugate.
Then $\Gamma$ is a discrete subgroup of $G$ because it is Zariski dense but
not analytically dense.  It is not analytically dense because the product
of two entries of $x$ and $\bar{x}$ in the same position is an ordinary
integer.  On the other hand the projection of $\Gamma$ onto either factor
is $\SL(2,\Z[\sqrt{3}])$, which is dense.  This is a counterexample not only
to \lem{l:closure}, but to its use in combination with \lem{l:diag} below.
\end{example}

\begin{lemma}[Diagonal; P. Hall \cite{Hall:eulerian}] Suppose that each of
$G_1, G_2, \ldots, G_\ell$ is a minimal simple Lie group or a non-abelian
finite simple group, and suppose that
$$H \subseteq G = G_1 \times G_2 \times \cdots \times G_\ell$$
is a closed subgroup that surjects onto each factor $G_k$.  Then $H$
is a diagonal subgroup of $G$.
\label{l:diag} \end{lemma}

A Lie group $G$ is \emph{minimal simple} if its Lie algebra $\mg$ is
simple, and if it is connected and has no center.  (See the remarks
after \lem{l:stem}.)  In the setting of the lemma, a subgroup $H$
is \emph{diagonal} if the following holds:  There exist isomorphisms
$\phi_{j,k}:G_j \to G_k$ for some pairs of the groups $G_1, G_2, \ldots,
G_\ell$.  For each $k$, at most one $\phi_{j,k}$ should be chosen, and
only when $j < k$.  Then $H$ consists of those elements $(g_1,g_2, \ldots,
g_\ell)$ such that $g_k$ is unrestricted when $\phi_{j,k}$ does not exist,
and such that $g_k = \phi_{j,k}(g_j)$ when $\phi_{j,k}$ does exist.

After showing that various representations $V = X(\vc)$ of $\Gamma = \PB_n$
are dense, \lem{l:diag} is the main tool to show that they are jointly
dense.  Although it is stated in a different form, this tool
is similar in spirit to the Independence Lemma of \cite{AAEL:tutte}.

\begin{lemma}[Stem] If $\pi:\tG \to G$ is a stem extension of groups,
and if $H \subseteq \tG$ projects onto all of $G$, then $H = \tG$.
\label{l:stem} \end{lemma}

One role of \lem{l:stem} is to clear up some confusion in the terminology
for Lie groups.  A surjective homomorphism $\pi:\tG \to G$ is a \emph{stem
extension} if $\ker \pi$ is central and lies in the commutator subgroup of
$\tG$.  A group is \emph{quasisimple} if it is a stem extension of a simple
group.  A Lie group $G$ is quasisimple if and only if it is connected and
its Lie algebra $\mg$ is simple.  However, the standard terminology for such
a Lie group is that it is ``simple".  Thus, $\SL(N,\C)$ is only quasisimple
as a group even though it is a simple Lie group.  We call a Lie group
\emph{minimal simple} if it is simple as a group, for example $\PSL(N,\C)$.

\begin{lemma}[Commutator denseness] Consider the following diagram
of groups:
$$\begin{tpic}[baseline=(current bounding box.east)]
\matrix (m) [matrix of math nodes,row sep=2em]
{\Gamma_1 & \subseteq & G_1 & \supseteq & G_3 & = [G_1,G_1]\\ 
    \Gamma_2 & \subseteq & G_2 & \supseteq & G_4 & = [G_2,G_2] \\ };
\draw[->>] (m-1-1) -- (m-2-1);
\draw[->>] (m-1-3) -- (m-2-3);
\draw[->>] (m-1-5) -- (m-2-5);
\end{tpic}.$$
Suppose that each $G_k$ is a topological group, that the projection $G_3
\twoheadrightarrow G_4$ is a stem extension that sends closed subgroups to
closed subgroups, that the groups on the right are commutator subgroups
as indicated, and that $\Gamma_2$ is the image of $\Gamma_1$ in $G_2$.
If the closure of $\Gamma_2$ contains $G_4$, then $\Gamma_1 \cap G_3$
is dense in $G_3$.
\label{l:comm} \end{lemma}

Note that the projection $G_1 \twoheadrightarrow G_2$ in \lem{l:comm} need
not be a stem extension.  The lemma will be used when the four groups $G_k$
are $\GL(N,\C)$, $\PGL(N,\C)$, $\SL(N,\C)$, and $\PSL(N,\C)$; or when they
are analogous real or unitary groups; or when they are products of these
groups.  It is useful in conjunction with \lem{l:diag} because for instance
$\PSL(N,\C)$ is a minimal simple Lie group, while $\SL(N,\C)$ is simple
but not minimal and $\GL(N,\C)$ is not simple.  Among other consequences,
we need only be interested in representations up to projective equivalence.

\begin{lemma}[Subrepresentation] Let $V$ be a representation of a group
$G$ over a field $\F$, and suppose that $V$ decomposes as a direct sum of
irreducible representations with multiplicity:
$$V = n_1 V_1 \oplus n_2 V_2 \oplus \cdots \oplus n_\ell V_\ell.$$
Then every $G$-invariant subspace $W \subseteq V$ decomposes
as
$$W \cong m_1 V_1 \oplus m_2 V_2 \oplus \cdots \oplus m_\ell V_\ell,$$
where for each $k$,
$$m_k V_k \cong W_k \subseteq n_k V_k.$$
\label{l:subrep} \end{lemma}

\lem{l:subrep} is used mainly for \lem{l:connect} below.

\begin{lemma}[Connectness] Suppose that a vector space $X$ over a field
$\F$ is a multiplicity-free representation of two groups $G$ and $H$.  Let
$$X = V_1 \oplus V_2 \oplus \cdots \oplus V_n$$
be the $G$-irreducible decomposition of $X$ and let 
$$X = W_1 \oplus W_2 \oplus \cdots \oplus W_\ell$$
be the $H$-irreducible decomposition.  Define a directed graph $C(X,G,H)$
on these summands, with an edge from $V_j$ to $W_k$ (or vice versa) if
there exists $v \in V_j$ with a non-zero component in $W_k$.  Then $X$
is an irreducible of the free product $G * H$ if and only if $C(X,G,H)$
is strongly connected.
\label{l:connect} \end{lemma}

To review, a directed graph $C$ is \emph{strongly connected} if for every
two vertices $p$ and $q$, there is a directed path in $C$ from $p$ to $q$.

We will use \lem{l:connect} in combination with \lem{l:adj} to inductively
show that the adjoint space $V = \sl(X(n \cdot 1,c,t))$ is irreducible under
$B_n$.  The groups $G$ and $H$ will be two copies of $B_{n-1}$ in $B_n$.

\begin{lemma}[Real Zariski] Let $G \subseteq \SL(N,\C)$ be a subgroup
which is real Zariski closed and complex Zariski dense.  Then either $G
= \SL(N,\C)$, or the connected subgroup of $G$ is a conjugate of either
$\SL(N,\R)$ or $\SU(P,Q)$ with $P+Q = N$, or $\SL(N/2,\H)$ when $N$ is even.
In particular, if the Lie algebra of $G$ contains $x$ and $ix$ for some $x$,
then $G = \SL(N,\C)$.
\label{l:real} \end{lemma}

The first claim of \lem{l:real} quickly implies cases 2 and 4 of \cor{c:real}
and real Zariski forms of cases 3 and 5.  The second claim of \lem{l:real}
can be viewed as a complex version of \lem{l:adj}.

\begin{lemma}[Rotation] Let $\Gamma$ be a subgroup of a Lie group $G$
that acts on a vector space $V$.  If $\Gamma$ has an elliptic element $g$
of infinite order, then $\Gamma$ is indiscrete in $G$ in the analytic
topology.  If $G = \SU(P,Q)$ with its defining representation and if
$\Gamma$ is finitely generated and analytically dense, then $\Gamma$
has an elliptic element $g$ of infinite order.
\label{l:rot} \end{lemma}

\lem{l:rot} is a tool to establish analytic denseness once the
real Zariski closure has been computed.

\begin{lemma}[Zassenhaus] If $\Gamma$ is a finitely generated discrete
group and $G$ is a semisimple real Lie group, then the set of dense
representations $\rho:\Gamma \to G$ (in the analytic topology) is an open
subset of the set of all representations.
\label{l:zassenhaus} \end{lemma}

We will not directly use \lem{l:zassenhaus} as a tool, but it sheds light
on what to expect in \thm{th:discrete}.

\section{Proofs of the tools and more lemmas}
\label{s:toproofs}

In this section, we will prove the lemmas in \sec{s:tools}.  We will also
discuss some other lemmas that are either related but not directly used,
or are used but are more technical or secondary.

\lem{l:adj} holds in both the analytic topology and in the Zariski
topology.  Its proof is elementary in light of the fact that a closed
subgroup of a Lie group is a Lie group.  The idea of stating the lemma
might not be considered elementary, although it is standard in some contexts.

\lem{l:zariski} is also elementary; it follows immediately from the
fact that a linear subspace of $V$ is Zariski closed, and that the group
action $G \times V \to V$ is Zariski continuous.  What is somewhat less
elementary is the context of other basic facts about the Zariski topology.
For instance, one relevant fact is that the group law $G \times G \to G$ is
Zariski continuous, \ie, $G$ is a topological group relative to its Zariski
topology.  This follows from the fact that matrix multiplication is Zariski
continuous, because there is a polynomial formula for matrix multiplication.

\lem{l:compact} is originally due to Chevalley and follows quickly
from the Stone-Weierstrass theorem \cite[Prop. 4.6.1]{Witte:ratner}.

\lem{l:closure} is a standard result in the theory of algebraic
groups.  It follows from another theorem of Chevalley, that the image of
any algebraic set under an algebraic map is constructible
\cite[Cor. 1.4]{Borel:gtm}.

\begin{proof}[Proof of \lem{l:diag}] The proof is by induction on $\ell$.
We first number the Lie factors of $G$ in order of non-increasing dimension,
and then after the Lie factors number the finite factors in order of
non-increasing cardinality. Let
$$\pi:H \to G' = G_1 \times G_2 \times \cdots \times G_{\ell-1}$$
be the projection of $H$ onto the first $\ell-1$ factors.  We can assume
by induction that $\pi(H) = G'$ by applying the lemma to $\pi(H)$ and then
replacing $G'$ by $\pi(H)$.

The kernel $\ker \pi$ is a subgroup of $G_\ell$, and we claim that it
is normal.  Suppose that $h = (1,1,\ldots,1,a) \in \ker \pi$.  Since $H$
surjects onto $G_\ell$, this means that for every $g_\ell \in G_\ell$,
there exists
$$g = (g_1,g_2,\ldots,g_\ell) \in H$$
for some choices of the other coordinates.  Then 
$$ghg^{-1} = (1,1,\ldots,g_\ell a g_\ell^{-1}) \in H,$$
so that $\ker \pi$ is normal in $G_\ell$ as claimed.  If $\ker \pi$ is
non-trivial, then $\ker \pi = G_\ell$ and $H = G$ and we are done.

If instead $\ker \pi = 0$, then $H$ is the graph of a continuous group
homomorphism $\alpha:G' \to G_\ell$, and $\alpha$ factors as a direct
product of homomorphisms $\alpha_k:G_k \to G_\ell$.  The domain and target
of $\alpha_k$ are both minimal simple Lie groups and either $\dim G_k \ge
\dim G_\ell$ or $|G_k| \ge |G_\ell|$.  Therefore $\alpha_k$ is either
the trivial homomorphism or an isomorphism.  Moreover, since $G_\ell$
is non-commutative, it is not possible for more than one $\alpha_k$ to be
surjective.  Thus at most one $\alpha_k$ is an isomorphism and the others
are the trivial homomorphism.  This gives $H$ the structure promised by
the lemma.
\end{proof}

\begin{proof}[Proof of \lem{l:stem}] Since $H$ projects onto $G$, it
intersects every coset of $K = \ker \pi$.  Therefore it is enough to show
that $H$ also contains $K$.  By hypothesis, there exist commutators $[a,b]
\in K$, with $a,b \in \tG$, that generate $K$.  Since $K$ is central,
$[a,b] = [ag,bh]$ for any $g,h \in K$.  Since $H$ meets every coset of
$K$, we can choose $g,h$ so that $ag, bh \in H$.  Thus $H$ contains $K$,
as desired.  \end{proof}

\begin{proof}[Proof of \lem{l:comm}] The main idea is 
to define the commutator map
$$c:G_2 \times G_2 \to G_4.$$
This map is continuous and surjective, so it sends the dense
set $\Gamma_2 \times \Gamma_2$ to dense set in $G_4$,
namely $\Gamma_4 = [\Gamma_2,\Gamma_2]$.  Then
$\Gamma_4$ is covered by $\Gamma_3 = [\Gamma_1,\Gamma_1]$.
By \lem{l:stem}, and by the hypothesis that the map 
from $G_3$ to $G_4$ sends closed subgroups to
closed subgroups, $\Gamma_3$ is dense in $G_3$.
At the same time, $\Gamma_3 \subseteq \Gamma_1$, as
desired.
\end{proof}

\lem{l:subrep} is a standard fact in representation theory.  We sketch
a proof for completeness, and because the lemma and proof are similar to
\lem{l:diag} and its proof.

\begin{proof}[Proof of \lem{l:subrep}] We first rename the decomposition
of $V$ as
$$V = V_1 \oplus V_2 \oplus \cdots \oplus V_\ell,$$
where the summands are irreducible but not necessarily inequivalent.

The proof is by induction on $\ell$.  Let
$$\pi:W \to V' = V_1 \oplus V_2 \oplus \cdots \oplus V_{\ell-1}$$
be the projection of $W$ onto the first $\ell-1$ summands.  We can assume
by induction that $\pi$ is surjective, by applying the lemma to $\pi(W)$
and replacing $V'$ by $\pi(W)$.  If $\ker \pi$ is non-trivial, it is a
submodule of $V_\ell$ and it therefore is $V_\ell$.  Thus if $\ker \pi$
is non-trivial, then $W = V$ and we are done.

Suppose instead that $\pi$ is a bijection between $W$ and $V'$.  Then $W$ is
the graph of a linear map $\alpha:V' \to V_\ell$.  The map $\alpha$ is the
direct sum of maps $\alpha_k:V_k \to V_\ell$.  By Schur's lemma, $\alpha_k
= 0$ when $V_k \not\cong V_\ell$.  (This part of Schur's lemma does not
require an algebraically closed field.)  Thus $\alpha$ is supported only
on those summands of $V'$ isomorphic to $V_\ell$, which is the structure
promised by the lemma.
\end{proof}

\begin{proof}[Proof of \lem{l:connect}] The idea of this lemma is that if
several assertions are related by a strongly connected set of inferences,
then the assertions are all equivalent.  For instance it is common to say
``$p$ implies $q$ implies $r$ implies $p$, therefore $p$, $q$, and $r$
are equivalent".

Let $Y$ be a subspace of $X$ which is invariant under both $A$ and $B$,
and consider the assertions that $V_j \subseteq Y$ or that $W_k \subseteq Y$.
When $C(X,G,H)$ has an edge from $V_j$ to $W_k$, then $V_j \subseteq Y$
implies $W_k \subseteq Y$, and vice-versa.  By hypothesis, these
implications are a strongly connected graph, so that $Y$ must
either contain all of the summand or none of them.  Thus $Y = X$
or $Y = \{0\}$, so that $X$ is irreducible.

On the other hand, if $C(X,G,H)$ is not strongly connected, then it has
a strongly connected component $D$ with no outward edges.  It is easy to
confirm that the summands $\{V_j\}$ of $D$ have the same direct sum as the
summands $\{W_k\}$ of $D$.  This direct sum $Y$ is then both $G$-invariant
and $H$-invariant, so that $X$ is not irreducible.
\end{proof}

\begin{proof}[Proof of \lem{l:real}] This lemma is almost just the
classification of real forms of the complex simple Lie algebra $\sl(N,\C)$,
but it needs some extra reasoning at the level of algebraic groups.  The real
forms are well-known and match the conclusion of the lemma; see for instance
Fulton and Harris \cite[\S 26.1]{FH:gtm} for a list at the Lie algebra level.

The real algebraic group $G$ has a complexification $G_\C$.  We can
assume, by passing to a finite-index subgroup, that $G$ is connected in
the real Zariski topology.  Then by construction, $G_\C$ is a subgroup
of the complexification $\SL(N,\C) \times \SL(N,\C)$ of $\SL(N,\C)$
as a real algebraic group.  $G_\C$ maps to each factor of $\SL(N,\C)$,
and each of these maps $\alpha$ and $\bar{\alpha}$ extends the inclusion
of $G$.  The map $\alpha$ is complex algebraic and $\bar{\alpha}$ is its
complex conjugate.  By \lem{l:closure}, the image of $\alpha$ is closed,
and then surjective by the hypothesis that $G$ is complex Zariski dense.
Thus by \lem{l:diag}, $G_\C$ is either all of $\SL(N,\C) \times \SL(N,\C)$,
or it is the graph of a real algebraic automorphism of $\SL(N,\C)$.
Thus either $G = \SL(N,\C)$, or $G_\C = \SL(N,\C)$ and $G$ is 
a real form, as desired.
\end{proof}

\begin{proof}[Proof of \lem{l:rot}]
The first claim of the lemma is elementary.  In a basis in which $g$ is
diagonal, it is contained in the compact group of diagonal unitary matrices.
Therefore if $g$ has infinite order, it generates an indiscrete subgroup.

The second claim rests on two ideas.  The first idea is to argue by counting
degrees of freedom that the elliptic elements in $\SU(P,Q)$ include
an open set.  Therefore any dense subgroup has many elliptic elements.
Let $N = P + Q$.  Then
$$\dim \SU(P,Q) = N^2 - 1.$$
Moreover, the set of orthogonal line bases in $\C^{(P,Q)}$, with $P$
positive lines and $Q$ negative lines, is a manifold of complex dimension
$\binom{N}2$ and real dimension $N^2-N$.  If $B$ is such a line basis, then
we can construct an elliptic element $g$ which is diagonal in this basis.
The element $g$ has $N-1$ additional degrees of freedom, because its entries
are complex numbers of norm 1 whose product is 1.  Indeed, for every $\eps >
0$, there is an open set of elliptic elements with at least one eigenvalue
$\exp(i \theta)$ with $\eps < \theta < 2\eps$.  Thus $\Gamma$ must have
elliptic elements of either infinite order or unbounded finite order.

The second idea is similar to the proof of Malcev's theorem that $\Gamma$
is residually finite.  The matrix entries of the elements of $\Gamma$
lie in a field $\F \subseteq \C$ which is finitely generated over $\Q$.
It is an interesting fact that any subfield of a finitely generated field
is still finitely generated \cite[Th. 3.3.5]{Roman:gtm}; in particular the
algebraic subfield of $\F$ is a finite-dimensional field $\F' \supseteq \Q$.
A root $\lambda$ of a characteristic polynomial of an element $g \in G$
is then either transcendental, or it lies in a field $\F'' \supseteq \F'$
whose degree over $\F'$ is bounded by $N$.  This imposes an upper bound
on the order of $\lambda$ if it is a root of unity, because for every
$n$, there are only finitely many roots of unity of algebraic degree $n$.
Thus it is not possible for $\Gamma$ to have elliptic elements of unbounded
finite order; it must instead have elliptic elements of infinite order.
\end{proof}

\lem{l:zassenhaus} follows quickly from the standard result
\cite[Th. 2.1]{BG:dense} that if $\Gamma$ has elements in a small enough
neighborhood of $1 \in G$ whose logarithms generate the Lie algebra $\mg$
(or even span $\mg$), then $\Gamma$ is dense.

\section{Quantum algebra}
\label{s:quantum}

In this section we review some of the properties of the Kauffman
bracket and by extension the Jones polynomial.  For a more
complete introduction to this theory, see for instance
Kauffman and Lins \cite{KL:recoupling}.

The Jones-Wenzl projector of color $c$ can be defined recursively 
as follows:
$$\begin{tpic}[style=basic]
\draw (0,-3) -- (0,3);
\draw[anchor=north] (0,-3.1) node {$c$};
\draw[thick,darkred,fill=white] (-2,-.25) rectangle (2,.25);
\end{tpic}
\;=\; \begin{tpic}[style=basic]
\draw (0,-3) -- (0,3);
\draw (-3,-3) -- (-3,3);
\draw[anchor=north] (0,-3.1) node {$c-1$};
\draw[thick,darkred,fill=white] (-2,-.25) rectangle (2,.25);
\end{tpic}
\;+ \frac{[c-1]}{[c]}\; \begin{tpic}[style=basic]
\draw (0,-3) -- (0,-2); \draw (0,2) -- (0,3);
\draw (1,-2) -- (1,2);
\draw (-1,-2) .. controls (-1,0) and (-3,0) .. (-3,-2) -- (-3,-3);
\draw (-1,2) .. controls (-1,0) and (-3,0) .. (-3,2) -- (-3,3);
\draw[anchor=north] (0,-3.1) node {$c-1$};
\draw[thick,darkred,fill=white] (-2,-2.25) rectangle (2,-1.75);
\draw[thick,darkred,fill=white] (-2,2.25) rectangle (2,1.75);
\draw[anchor=west] (1.1,0) node {$c-2$};
\end{tpic}\;.
$$
Here and below, a strand labelled by $n$ is shorthand for $n$ strands.  Also
$$[n] = \frac{t^{n/2} - t^{-n/2}}{t^{1/2}-t^{-1/2}}$$
is a \emph{quantum integer}.  Because of a quantum integer appearance in the
denominator, if $t$ has order $r$, then
the Jones-Wenzl projector of color $c$ is only defined when $c \le r-1$.
Moreover, the color $r-1$ is suppressed, because any reduced
skein space such as $X(n \cdot 1,r-1,t)$ vanishes.  So we say that $c$
is an \emph{admissible color} if $c \le r-2$.

If $\vc = (c_1,c_2,\ldots,c_\ell)$ is a vector of admissible colors,
then we can define a generalized skein space $W(\vc,t)$ and
reduced skein space $X(\vc,t)$.  The space $W(\vc,t)$ is defined
using $\ell$ clasps of the respective colors:
$$\begin{tpic}[basic,scale=4]
\draw (1.28,0.24) .. controls +(130.6:.5) and +(236.3:.5) .. (1.4,2.2);
\draw (1.56,0.48) .. controls +(130.6:.5) and +(236.3:.5) .. (1.6,2.05);
\draw (1.84,0.72) .. controls +(130.6:.5) and +(236.3:.5) .. (1.8,1.9);
\draw (2.12,0.96) .. controls +(130.6:.3) and +(236.3:.3) .. (2,1.75);
\draw (1.2,2.35) .. controls +(236:.5) and +(31:.5) .. (.45,.45);
\draw (.4,2.2) .. controls +(-45:.5) and +(31:.5) .. (.3,.7);
\draw (.2,2.0) .. controls +(-45:.5) and +(31:.5) ..  (.15,.95);
\draw[thick,darkred] (0,1.2) -- (.6,.2);
\draw[thick,darkred] (.1,1.9) -- (.5,2.3);
\draw[thick,darkred] (1,2.5) -- (2.2,1.6);
\draw[thick,darkred] (1,0) -- (2.4,1.2);
\end{tpic}$$
As before, the reduced space $X(\vc,t)$ is defined by gluing together
two skeins with projectors and dividing by the kernel of the resulting
bilinear form.

We will also need three standard results in the skein theory of the
Kauffman bracket and one structural result.  We omit the proofs; see
\cite{KL:recoupling} for the general theory.

The structural result that we will need is the following well-known
splitting formula.

\begin{theorem}[Splitting]
If $\vc$ and $\vc'$ are two lists of colors and $\vc \oplus \vc'$
is their concatenation, then
\eq{e:gluing}{X(\vc \oplus \vc',t) \cong \bigoplus_a X(\vc,a,t)
    \tensor X(a,\vc',t),}
where the map
$$X(\vc,a,t) \tensor X(a,\vc',t) \to X(\vc \oplus \vc',t)$$
is given by gluing together skeins using a projector of weight $a$.
\label{th:split} \end{theorem}

We will need the effect of a full twist on a projector. 

\begin{lemma}
$$\begin{tpic}[style=basic]
\draw[style=tangle] (3,1) .. controls (3,-1) and (0,0) .. (0,4);
\draw[style=tangle] (0,-2) .. controls (0,2) and (3,3) .. (3,1);
\draw (0,-4) -- (0,-2);
\draw (0,-4.1) node[anchor=north] {$c$};
\draw[thick,darkred,fill=white] (-1.5,-2.25) rectangle (1.5,-1.75);
\end{tpic}\; = t^{c(c+2)/4}\;\begin{tpic}[style=basic]
\draw (0,-4) -- (0,4);
\draw[anchor=north] (0,-4.1) node {$c$};
\draw[thick,darkred,fill=white] (-1.5,-.25) rectangle (1.5,.25);
\end{tpic}$$
\label{l:twist} \end{lemma}

We will also want an explicit change-of-basis formula between the two ways
to split $X(c,1,c,1,t)$.

\begin{lemma} Let
\begin{align*}
v_1 &\;=\; \begin{tpic}[style=basic]
\draw (-.5,-3) .. controls (-.5,-1) and (.5,1) .. (.5,3);
\draw (-3,0) .. controls (-1.5,0) and (-.5,1) .. (-.5,3);
\draw (3,0) .. controls (1.5,0) and (.5,-1) .. (.5,-3);
\draw[thick,darkred] (-1.5,3) -- (1.5,3);
\draw[thick,darkred] (-1.5,-3) -- (1.5,-3);
\draw (-3,0) node[anchor=east] {$1$}; \draw (3,0) node[anchor=west] {$1$};
\draw (0,3) node[anchor=south] {$c$}; \draw (0,-3) node[anchor=north] {$c$};
\begin{scope}[rotate=-40]
\draw[thick,darkred,fill=white] (-2.5,-.25) rectangle (2.5,.25);
\end{scope} \end{tpic} &
v_2 &\;=\; \begin{tpic}[style=basic]
\draw (-.5,3) .. controls (-.5,1) and (.5,-1) .. (.5,-3);
\draw (-3,0) .. controls (-1.5,0) and (-.5,-1) .. (-.5,-3);
\draw (3,0) .. controls (1.5,0) and (.5,1) .. (.5,3);
\draw[thick,darkred] (-1.5,3) -- (1.5,3); 
\draw[thick,darkred] (-1.5,-3) -- (1.5,-3);
\draw (-3,0) node[anchor=east] {$1$}; \draw (3,0) node[anchor=west] {$1$};
\draw (0,3) node[anchor=south] {$c$}; \draw (0,-3) node[anchor=north] {$c$};
\begin{scope}[rotate=-40]
\draw[thick,darkred,fill=white] (-2.5,-.25) rectangle (2.5,.25);
\end{scope} \end{tpic} \\
w_1 &\;=\; \begin{tpic}[style=basic]
\draw (-.5,-3) .. controls (-.5,-1) and (.5,1) .. (.5,3);
\draw (-3,0) .. controls (-1.5,0) and (-.5,1) .. (-.5,3);
\draw (3,0) .. controls (1.5,0) and (.5,-1) .. (.5,-3);
\draw[thick,darkred] (-1.5,3) -- (1.5,3);
\draw[thick,darkred] (-1.5,-3) -- (1.5,-3);
\draw (-3,0) node[anchor=east] {$1$}; \draw (3,0) node[anchor=west] {$1$};
\draw (0,3) node[anchor=south] {$c$}; \draw (0,-3) node[anchor=north] {$c$};
\begin{scope}[rotate=40]
\draw[thick,darkred,fill=white] (-2.5,-.25) rectangle (2.5,.25);
\end{scope} \end{tpic}&
w_2 &\;=\; \begin{tpic}[style=basic]
\draw (-.5,3) .. controls (-.5,1) and (.5,-1) .. (.5,-3);
\draw (-3,0) .. controls (-1.5,0) and (-.5,-1) .. (-.5,-3);
\draw (3,0) .. controls (1.5,0) and (.5,1) .. (.5,3);
\draw[thick,darkred] (-1.5,3) -- (1.5,3); 
\draw[thick,darkred] (-1.5,-3) -- (1.5,-3);
\draw (-3,0) node[anchor=east] {$1$}; \draw (3,0) node[anchor=west] {$1$};
\draw (0,3) node[anchor=south] {$c$}; \draw (0,-3) node[anchor=north] {$c$};
\begin{scope}[rotate=40]
\draw[thick,darkred,fill=white] (-2.5,-.25) rectangle (2.5,.25);
\end{scope} \end{tpic}.
\end{align*}
Then
\begin{align*}
v_1 &= \frac{[c][c+2]}{[c+1]^2} w_1 + \frac{1}{[c+1]} w_2 \\
v_2 &= -\frac{1}{[c+1]} w_1 + w_2.
\end{align*}
\label{l:cob} \end{lemma}

Finally we will need the following two dimension inequalities:

\begin{lemma} Let $t$ be a root of unity of order $r$
with $r$ even, and let $0 \le c < \frac{r-2}2$.   
If $\dim X(n \cdot 1,c,t) > 0$, then
$$\dim X(n \cdot 1,c,t) > \dim X(n \cdot 1,r-2-c,t).$$
\label{l:dim1} \end{lemma}

\begin{proof} \thm{th:split} implies a recursive characterization
of the numbers
$$d(n,c,r) = \dim X(n \cdot 1,c,t)$$
when $c$ is admissible. Namely,
\begin{align*}
d(0,0,r) &= 1 \\
d(0,c,r) &= 0 \qquad \qquad (0 < c) \\
d(n,c,r) &= d(n-1,c-1,r) + d(n-1,c+1,r) \\
& \qquad \qquad \qquad (0 < c < r-2) \\
d(n,0,r) &= d(n-1,1,r) \\
d(n,r-2,r) &= d(n-1,r-3,r)
\end{align*}
A simple induction argument shows that $d(n,c,r) > 0$ when $c \le n$
and $c+n$ is even and that $d(n,c,r) = 0$ otherwise.  Now let
$$e(n,c,r) = d(n,c,r) - d(n,r-2-c,r).$$
It is easy to check that $e(n,c,r)$ satisfies the same recurrence as
$d(n,c,\frac{r}2)$.  Therefore these numbers are equal, and the inequality
for $d(n,c,r)$ is the desired claim.
\end{proof}

\begin{lemma} If $1 \le \dim X(n \cdot 1,c,t) \le 2$ and $c \ge 5$,
then either $c = n$ or $n \le 3$ or $(n,c) = (0,4)$.
\label{l:dim2} \end{lemma}

\begin{proof} We use the abbreviation $d(n,c,r)$ in \lem{l:dim1}
and the recurrence that these numbers satisfy. We obtain 
\begin{align*}
d(4,0,r) &= 2 & d(4,2,r) &= 3 \\
d(5,1,r) &= 5 & 3 \le d(5,3,r) &\le 4.
\end{align*}
It is then easy to check the lemma by induction for $n \ge 6$.
\end{proof}

\section{Proof of the main results}
\label{s:mproofs}

The proof of \thm{th:complex} is by mutual induction with \thm{th:joint}.
For convenience, we let $r$ be the order of $t$ if it is a root of unity,
and let $r = \infty$ otherwise.

\begin{proof}[Proof of \thm{th:joint}] \lem{l:comm} tells us that we can
replace each group $G_k$ with $\PSL(X)$, $\PSL(X_\R)$, or $\PSU(X)$, where
$X = X(n \cdot 1,c,t)$.  \lem{l:closure} tells us that we can replace the
action of $B_n$ by its Zariski closure in each factor.  Then \lem{l:diag}
tells us that our aim is to show that two representations $X = X(n \cdot
1,c,t)$ and $X' = X(n \cdot 1,c',t)$ are inequivalent, in the sense that
there does not exist an isomorphism $\alpha$ of the corresponding groups
that equates the actions of $B_n$.  Moreover, the graph of $\alpha$
has to be closed in the relevant Zariski topology, which implies
that $\alpha$ is continuous in the analytic topology.

Consider first the case
$$\alpha:\PSL(X) \to \PSL(X').$$
In this case, the only associated automorphisms of $\PGL(N,\C)$ are:
\begin{align*}
\alpha(x) &= x  & \alpha(x) &= (x^{-1})^T \\
\alpha(x) &= \bar{x} & \alpha(x) &= (\bar{x}^{-1})^T
\end{align*}
Here $x$ is a matrix.  If instead $\alpha$ is a map from $\PSU(X)$ to
$\PSU(X')$, or from $\PSL(X_\R)$ to $\PSL(X'_\R)$, then the automorphisms
have the same formulas, except that the last two are equivalent to the
first two.  In the complex Zariski topology, only the first two choices of
$\alpha$ have a closed graph.  Finally, since we are working in projective
linear groups, $x$ is only defined up to a scalar factor.

So we are interested in a ``fingerprint" of the representation $X$ that
will distinguish it from other representations $X'$.  the fingerprint
can make use of the spectrum of the action of an element $g \in B_n$,
up to a scalar factor, up to inversion, and up to conjugation, because
that information is preserved by all choices of $\alpha$.  If $\dim X =
1$, then fingerprinting is not necessary because the corresponding group
such as $\PSL(X)$ is trivial.

By \thm{th:split},
$$X(n \cdot 1,c,t) \cong X((n-1) \cdot 1,c+1,t)
    \oplus X((n+1) \cdot 1,c-1,t).$$
Let $g$ be the full twist on the first $n-1$ strands.  Then by \lem{l:twist},
its eigenvalues are proportional to $t^{(c+1)(c+3)/4}$ on $X((n-1) \cdot
1,c+1,t)$ and $t^{(c-1)(c+1)/4}$ on $X((n-1) \cdot 1,c-1,t)$.  Or, after
rescaling, the eigenvalues are $t^{c+1}$ and $1$.  One of the eigenvalues
is suppressed when $c = 0$ and when $c = r-2$.  Together with $n$ itself,
this is a complete fingerprint for $X$ when $t$ is not a root of unity
or when $r$ is odd.  When $r$ is even, this data does not distinguish
$c$ from $r-2-c$, because we can switch the spectrum $(t^{c+1},1)$ with
$(t^{r-2-c+1},1)$ by inverting $t$ and rescaling.  However, by \lem{l:dim1},
$$\dim X(n \cdot 1,c,t) > \dim X(n \cdot 1,r-2-c,t).$$
Thus the strand number $n$, the eigenvalues of $g$, and the dimension
$\dim X$ are a complete fingerprint for $X$ and prevent the existence of
the isomorphism $\alpha$.
\end{proof}

\begin{proof}[Proof of \thm{th:complex}] The proof is by induction and
the base case is $X(3 \cdot 1,1,t) = X(4 \cdot 1,0,t)$.  First suppose
that $t$ is not a root of unity.  In this case, let $\tau_1$ and $\tau_2$
be the braid generators of $B_3$ acting on $X(3 \cdot 1,1,t)$.   Then the
action of $\tau_1$ is
$$ \begin{pmatrix} t^{-3/4} & 0 \\ 0 & -t^{1/4} \end{pmatrix} $$
in the basis
$$\begin{tpic}[style=basic,scale=1.5]
\draw (0,0) .. controls (2,0) and (2,1) .. (0,1);
\draw (0,2) .. controls (2,2) and (2,3) .. (0,3);
\draw[darkgray,dashed] (0,3.5) -- (0,-.5);
\end{tpic} \qquad \qquad
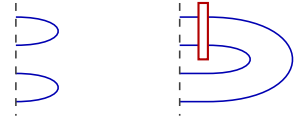
\begin{tpic}[style=basic,scale=1.5]
\draw (-1,0) -- (0,0) .. controls (4,0) and (4,3) .. (0,3) -- (-1,3);
\draw (-1,1) -- (0,1) .. controls (2,1) and (2,2) .. (0,2) -- (-1,2);
\draw[darkgray,dashed] (-1,3.5) -- (-1,-.5);
\draw[thick,darkred,fill=white] (-.333,1.5) rectangle (0,3.5);
\end{tpic}$$
while the action of $\tau_2$ has the same matrix as $\tau_1$,
but in the basis
$$\begin{tpic}[style=basic,scale=1.5]
\draw (0,0) .. controls (4,0) and (4,3) .. (0,3);
\draw (0,1) .. controls (2,1) and (2,2) .. (0,2);
\draw[darkgray,dashed] (0,3.5) -- (0,-.5);
\end{tpic} \qquad \qquad
\begin{tpic}[style=basic,scale=1.5]
\draw (-1,0) -- (0,0) .. controls (2,0) and (2,1) .. (0,1) -- (-1,1);
\draw (-1,2) -- (0,2) .. controls (2,2) and (2,3) .. (0,3) -- (-1,3);
\draw[darkgray,dashed] (-1,3.5) -- (-1,-.5);
\draw[thick,darkred,fill=white] (-.333,.5) rectangle (0,2.5);
\end{tpic}$$
By \lem{l:cob}, the change of basis matrix is:
$$M = \begin{pmatrix} -\frac{1}{[2]} & 1 \\ 
    \frac{[3]}{[2]^2} & \frac{1}{[2]} \end{pmatrix}.$$
The Zariski closure in $\PSL(2,\C)$ of the action of $\tau_1$, in its basis, is
$$\left\{ \begin{pmatrix} z & 0 \\ 0 & z^{-1} \end{pmatrix}
    \right\}_{z \in \C \setminus \{0\}}.$$
The Lie algebra $\sl(2,\C)$ has a multiplicity-free decomposition under this
action, and by \lem{l:subrep} it has only four $\tau_1$-invariant subspaces:
$$\begin{pmatrix} * & 0 \\ 0 & * \end{pmatrix} \qquad
\begin{pmatrix} * & * \\ 0 & * \end{pmatrix} \qquad
\begin{pmatrix} * & 0 \\ * & * \end{pmatrix} \qquad
\begin{pmatrix} * & * \\ * & * \end{pmatrix}.
$$
Because the change-of-basis matrix $M$ is full, the Lie algebra of the
Zariski closure of $\tau_2$ only lies in the last choice.  (This argument
is actually a special case of \lem{l:connect}.  Thus the Lie algebra of
the Zariski closure is all of $\sl(2,\C)$, and the group is $\PSL(2,\C)$.

Suppose that $t$ is a root of unity.  For convenience we instead
consider the action of $B_4$ on $X = X(4 \cdot 1,0,t)$.  Also for convenience,
we assume that $t$ is a principal root of unity, which is not a loss
of generality because the complex Zariski topology is Galois-invariant.
Let $\tau_1, \tau_2, \tau_3$ be the braid generators, and let
$$\sigma_3 = \tau_2 \tau_1 \qquad \sigma_4 = \tau_3 \tau_2 \tau_1.$$
Then in $\PU(X) \cong \PSU(X) \cong \SO(3)$, $\sigma_3$ has
order 3 and $\sigma_4$ has order 2.  Meanwhile
$$\tau_3 = \sigma_4 \sigma_3^{-1}$$
has eigenvalue ratio $-t$, so it has order $s = 2r$ when $r$ is odd, $s =
r/2$ when $r = 4k+2$, and $s = r$ when $r = 4k$.  The braids $\sigma_3$,
$\sigma_4$, and $\tau_3$ generate $B_4$, and their action, if finite,
is that of the $(2,3,s)$ triangle group.  However, this triangle group is
infinite when $s \ge 7$, so the action is dense.

When $r = 10$, then $s = 5$, and the action of $B_4$ is that of the
icosahedral group.   Crucially, this action is adjoint-irreducible and
the image is a simple group.  These properties will make it usable as a
base for induction even though it is not a case of the theorem.
The inductive step of this case will be saved for last.

For the inductive step, let $X = X(n \cdot 1,c,t)$.  If $c = n$, there is
nothing to prove because $\dim X = 1$.  If $c = 0$ and $n > 0$, there is
also nothing to do because
$$X(n \cdot 1,0,t) = X((n-1) \cdot 1,1,t).$$
Likewise if $t$ has order $r$ and $c = r-2$, then
$$X(n \cdot 1,c,t) = X((n-1) \cdot 1,c-1,t).$$
So suppose that $n > c > 0$ and that either $t$ is a non-root-of-unity or that
$c \le r-3$.  Then as in the proof of \thm{th:joint}.
$$X(n \cdot 1,c,t) \cong X((n-1) \cdot 1,c+1,t) \oplus X((n+1) \cdot 1,c-1,t).$$
Moreover, there are two different splittings, depending on whether we
restrict to the first $n-1$ strands or the last $n-1$ strands.  Let $G$
be the braid group on the first $n-1$ strands and let $H$ be the braid
group on the last $n-1$ strands.  Let $V_1$ and $V_2$ be the $G$-invariant
summands with colors $c+1$ and $c-1$, and likewise let $W_1$ and $W_2$
be the $H$-invariant summands.

We want to determine the position of $V_1$ and $V_2$ relative to $W_1$
and $W_2$ using the splitting
$$X \cong \bigoplus_a X((n-2) \cdot 1,a,t) \tensor X(a,1,c,1),$$
where the first factor uses the middle $n-2$ strands.  In particular
we will use the summand with $a = c$, which must be non-zero given our
assumptions on $n$.  Call this summand $Y \tensor Z$, where
$$Y = X((n-2) \cdot 1,c,t) \qquad Z = X(c,1,c,1).$$
Then the four skeins in \lem{l:cob} are in $Z$; when tensored with $Y$, they
place copies of $Y$ in each of $V_1, V_2, W_1, W_2$.  \lem{l:cob} then says
that the change-of-basis matrix between these skeins has no zero entries.

Thus the graph $C(X,G,H)$ of \lem{l:connect} is strongly connected, so
that $X$ is an irrep of $B_n$.  However, we are interested in $\sl(X)$,
which is a more complicated case.  This vector space decomposes as:
\eq{e:slx}{\sl(X) \cong \sl(V_1) \oplus \sl(V_2) \oplus (V_1 \tensor V_2^*)
    \oplus (V_2 \tensor V_1^*) \oplus I_G,}
where the last summand $I_G$ correspond to traceless operators on $X$
that act by scalars on both $V_1$ and $V_2$.  By induction, and by
\lem{l:zariski}, each term in this decomposition is $G$-irreducible, but
we would also like to know that it is multiplicity-free.  For that purpose,
\thm{th:joint} says that we can use the action of $\SL(V_1) \times \SL(V_2)$
to distinguish the summands.  Note also that if $\dim V_k = 1$,
then $\sl(V_k)$ vanishes; it is a null term in the decomposition.

Recall that the adjoint representation of $\SL(V)$ is irreducible
and inequivalent to the defining representation, and that the defining
representation is not self-dual if $\dim V \ge 3$.  Thus all of the terms
in the decomposition are inequivalent, except that $V_1 \tensor V_2^*$ and
$V_2 \tensor V_1^*$ are equivalent when $V_1$ and $V_2$ have dimension at
most 2.  By \lem{l:dim2}, the exception occurs only when $n = 4$ and $c = 2$.

To address this exception, let $g \in G$ be the full twist on the first $n-1$
strands, as in the proof of \thm{th:joint}.  The eigenvalue ratio of $g$
acting on $V_1$ and $V_2$ is $t^{c+1} = t^3$.  Therefore its eigenvalue
ratio on $V_1 \tensor V_2^*$ and $V_2 \tensor V_1^*$ is $t^6$.  We know
that $t^6 \ne 1$ since $t$ is a non-lattice root of unity, so $V_1 \tensor
V_2^*$ and $V_2 \tensor V_1^*$ are distinguished by the action of $g$.
Thus equation~\eqref{e:slx} is a multiplicity-free decomposition in all
germane cases.

By the same reasoning,
$$\sl(X) \cong \sl(V_1) \oplus \sl(V_2) \oplus (V_1 \tensor V_2^*)
    \oplus (V_2 \tensor V_1^*) \oplus I_H$$
is the $H$-irreducible decomposition of $\sl(X)$ and it is also
multiplicity-free.  Thus \lem{l:connect} applies to $\sl(X)$,
provided that $C(\sl(X),G,H)$ is strongly connected.  This is our
final claim to establish the theorem.

Let $v_1 \tensor v_2^* \in V_1 \tensor V_2^*$ be a vector, and 
consider its $H$-decomposition
$$v_1 \tensor v_2^* = w_1 \tensor w_1^* +
    w_1 \tensor w_2^* + w_2 \tensor w_1^* + w_2 \tensor w_2^*.$$
By the structure of $C(X,G,H)$, the vector $v_1$ and the dual vector $v_2^*$
can be chosen so that all four terms in the $H$-decomposition are non-zero.
Moreover, $w_k \tensor w_k^*$ must have a non-zero component in $\sl(V_k)$
when $\dim V_k \ge 2$, because it is a rank 1 operator and cannot be
proportional to the identity.  Thus in $C(\sl(X),G,H)$, there is an edge from
$V_1 \tensor V_2^*$ to every $H$-invariant summand other than possibly $I_H$.

We also claim that there is an edge from $I_G$ to at least one
$H$-irreducible summand other than $I_H$.  This will happen unless $I_G
= I_H$.  They cannot be equal, because the eigenspaces of $x \in I_G$ are
$V_1$ and $V_2$, while the eigenspaces of $y \in I_H$ are $W_1$ and $W_2$;
and these subspaces of $X$ are different.

Moreover, $\sl(X)$ is a self-dual representation using the bilinear
form $\Tr(xy)$ (the Killing form).  This implies that if $C(\sl(X),G,H)$
has an edge from a vertex $A$ to a vertex $B$, it also has an edge from
$B^*$ to $A^*$.  Moreover, for all of the edges constructed so far, we
can switch $V_1$ and $V_2$ with $W_1$ and $W_2$.  All told, these edges
render $C(\sl(X),G,H)$ strongly connected.  This concludes the proof
when $r \ne 10$.

\begin{figure}[htb]
\begin{scriptsize}
\begin{verbatim}
# Set up rings and variables.

R.<s> = QQ[]
RX.<x> = PolynomialRing(Frac(R))
t = s^4

# Conveniences

def mat(a): return matrix(Frac(R),a).transpose()
def comm(a,b): return a*b*a.inverse()*b.inverse()
    
# These are the generators of the Jones representation
# of B_4 acting on the skein space 1,1,1,1,2.

tau1 = mat([[s^(-3),0,0],[-s^(-1),-s,0],[0,0,-s]])
tau2 = mat([[-s,-s^(-1),0],[0,s^(-3),0],[0,-s^(-1),-s]])
tau3 = mat([[-s,0,0],[0,-s,-s^(-1)],[0,0,s^(-3)]])

# Make the commutator [tau2,tau2*tau3^3*tau2*tau1^(-1)].
# Verify its characteristic polynomial; it should print 0

c = comm(tau2,tau2*tau3^3*tau2*tau1^(-1))
print c.charpoly() - (x-1)*(x^2+(t-1+t^(-1))^3*x + 1)
\end{verbatim}
\end{scriptsize}
\caption{Sage code to compute $[\tau_2,\tau_2\tau_3^3\tau_2\tau_1^{-1}]$}
\label{f:sage} \end{figure}

Finally let $r = 10$.  The projective action of $B_3$ on $X(3 \cdot 1,1,t)$
is that of the icosahedral group.  Since this action is both
adjoint-irreducible and a finite simple group, it can be used
in \thm{th:joint}, and in the rest of the above argument
when $n = 4$.  We can conclude that the action of $B_4$
on $X(4 \cdot 1,2,t)$ is adjoint-irreducible, but we also claim
that it is indiscrete, and this claim needs a separate argument. An ad
hoc search in Sage reveals that commutator
$$g = [\tau_2,\tau_2\tau_3^3\tau_2\tau_1^{-1}]$$
has characteristic polynomial
$$\chi_g(x) = (x-1)(x^2+(t-1+t^{-1})^3 x + 1)$$
in its action on $X(4 \cdot 1,2,t)$, for any $t$ for which this
space is 3-dimensional.  (See \fig{f:sage}.)  When $t = \exp(\pi i/5)$,
the roots $\lambda$ other than $1$ have the form $\exp(i \theta)$
with 
$$|\theta| \approx 96.778652^\circ.$$
We claim that this is an irrational angle, although at first glance it is
not clear.  By Galois theory, $\lambda$ has degree at most 4 over $\Q$.
Thus if $\lambda$ did have finite order $s$, the Euler totient $\phi(s)$ of
$s$ would be at most 4, so that $s \le 10$.  The approximation to $\theta$
thus tells us that $g$ has infinite order and the action of $B_4$ is dense.

The projective action of $B_4$ on $X(4 \cdot 1,0,t)$ is also that of the
icosahedral group.  However, this and the action on $X(4 \cdot 1,2,t)$
are enough to show adjoint irreducibility when $n = 5$.  When $n = 5$
none of the actions can be discrete, and the case $r = 10$ thus merges
with the other cases.
\end{proof}

\begin{proof}[Proof of \cor{c:real}] As suggested in \sec{s:tools}, cases
2 and 4 follow from \lem{l:real} together with the existence of $X(n \cdot
1,c,t)_\R$ when $t$ is real, and the invariant Hermitian structure on $X(n
\cdot 1,c,t)$ when $|t| = 1$ as described in \sec{s:quantum}.  \lem{l:real}
likewise implies cases 3 and 5 follow in the real Zariski topology.

To complete cases 3 and 5, note that the criterion of \lem{l:rot}
is easily satisfied, albeit for different reasons.  In case 3, the ratio
of the eigenvalues of a braid generator $\tau_1$ (say) is
$-t$.  Thus when $|t| = 1$ but $t$ is not a root of unity, then
$\tau_1$ is an elliptic element of infinite order.  In case 5,
if $B_n$ acts discretely, then it acts finitely and is not
complex Zariski dense in $\SL(X(n \cdot 1,c,t))$.  As mentioned,
case 5 also follows from \lem{l:compact}.

As mentioned, case 6 follow from \lem{l:adj}. 

Finally consider case 1.  This case is established by induction on $n$.
First, let $X = X(3 \cdot 1,1,t)$.  By \lem{l:real}, we wish to show that
the projective action of $B_3$ is not contained in any copy of $\PSU(2)$
or $\PSU(1,1)$ or $\PSL(2,\R)$.  Specially to dimension 2, $\PSU(1,1)$ and
$\PSL(2,\R)$ are conjugate in $\PSL(2,\C)$.  Moreover, the normalizer of
$\PSL(2,\R)$ can be written as $\PGL(2,\R)$ in $\PGL(2,\C) = \PSL(2,\C)$.
The normalizer of $\PSU(2)$ is itself.

Again, the eigenvalue ratio of the braid generator $\tau_1$ is $-t$.
Since the ratio does not have norm 1, it is not contained in any conjugate
of $\PSU(2)$.  Since the ratio is not real, it is not contained in any
conjugate of $\PGL(2,\C)$.

In the inductive case, let $X = X(n \cdot 1,c,t)$.  This $X$ is a direct sum
of many copies of $X(3 \cdot 1,1,t)$ and $X(3 \cdot 1,3,t)$.  The commutator
subgroup $[B_3,B_3]$ of $B_3$ must still be real Zariski dense in its
action on $X(3 \cdot 1,1,t)$, while its action on $X(3 \cdot 1,3,t)$
is trivial.  Therefore the Lie algebra of the real Zariski closure of the
action of $[B_3,B_3]$ includes both $x$ and $ix$ for some $x \in \sl(X)$.
By \lem{l:real}, this implies that the action of $B_n$ is real Zariski dense.
\end{proof}

\section{Other results}
\label{s:other}

\begin{proof}[Proof of \thm{th:discrete}] It will be more convenient to
consider the action of $B_4$ on $X(4 \cdot 1,0,t)$.  Case 2 is a special
case of \thm{th:complex}, given that $\SU(X(4 \cdot 1,0,t))$ is compact.

\begin{figure*}[htb]
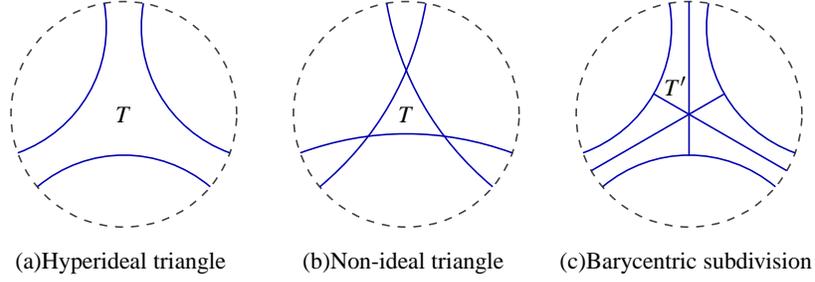

\begin{center}
\subfigure[Hyperideal triangle]{
\begin{tpic}[style=basic,scale=6]
\useasboundingbox (-1.2,-1.2) rectangle (1.2,1.2);
\draw[darkgray,dashed] (0,0) circle (1);
\draw (100:1) arc (10:-70:1.19175); \draw (220:1) arc (130:50:1.19175);
\draw (-20:1) arc (250:170:1.19175); \draw (0,0) node {$T$};
\end{tpic}}
\subfigure[Non-ideal triangle]{
\begin{tpic}[style=basic,scale=6]
\useasboundingbox (-1.2,-1.2) rectangle (1.2,1.2);
\draw[darkgray,dashed] (0,0) circle (1);
\draw (80:1) arc (-10:-50:2.7475); \draw (200:1) arc (110:70:2.7475);
\draw (-40:1) arc (230:190:2.7475); \draw (0,0) node {$T$};
\end{tpic}}
\subfigure[Barycentric subdivision]{
\begin{tpic}[style=basic,scale=6]
\useasboundingbox (-1.2,-1.2) rectangle (1.2,1.2);
\draw[darkgray,dashed] (0,0) circle (1);
\draw (100:1) arc (10:-70:1.1917); \draw (220:1) arc (130:50:1.1917);
\draw (-20:1) arc (250:170:1.1917); \draw (-.12,.25) node {$T'$};
\draw (270:0.364) -- (90:1); \draw (150:0.364) -- (-30:1);
\draw (30:0.364) -- (210:1);
\end{tpic}}
\end{center}
\caption{Hyperbolic triangles}
\label{f:triangle} \end{figure*}

To analyze cases 1 and 3, it is useful to work in $$\PGL(2,\R) \subseteq
\PGL(2,\C) = \PSL(2,\C),$$ and to recall that $\PGL(2,\R)$ is the isometry
group of the hyperbolic plane $\H^2$.  (Note that in \lem{l:real}, ``$\H$''
refers instead to the quaternions.)  The non-trivial elements with positive
determinant are rotations, hyperbolic translations, and parabolic motions.
The elements with negative determinant are reflections and glide reflections.
Recall also that $\PSL(2,\C)$ is the rotation group of hyperbolic space
$\H^3$.  Reflections in $\PGL(2,\R)$ are realized as rotations in $\H^3$ that
flip over $\H^2$.  In particular, reflections are conjugate in $\PSL(2,\C)$
to rotations by $\pi$, even though they are not conjugate in $\PGL(2,\R)$.

An element $g$ can be analyzed in terms of its eigenvalue ratio $\rho =
\lambda_1/\lambda_2$.  If $\rho > 0$, then $g$ is hyperbolic.  If $\rho =
1$ (formally, using generalized eigenvalues), then $g$ is the identity
or it is parabolic.  if $\rho = \exp(i\theta)$, then $g$ is a rotation
by an angle of $\theta$.  If $\rho = -1$, then $g$ is a reflection, and
otherwise if $\rho < 0$ then $g$ is a glide reflection.

As before, we let $\tau_1$, $\tau_2$, and $\tau_3$ be the braid
generators, and we let
$$\sigma_3 = \tau_2 \tau_1 \qquad \sigma_4 = \tau_3 \tau_2 \tau_1.$$
In the basis
$$\begin{tpic}[style=basic,scale=1.5]
\draw (0,0) .. controls (4,0) and (4,3) .. (0,3);
\draw (0,1) .. controls (2,1) and (2,2) .. (0,2);
\draw[darkgray,dashed] (0,3.5) -- (0,-.5);
\end{tpic}
\qquad \qquad
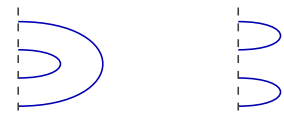
\begin{tpic}[style=basic,scale=1.5]
\draw (0,0) .. controls (2,0) and (2,1) .. (0,1);
\draw (0,2) .. controls (2,2) and (2,3) .. (0,3);
\draw[darkgray,dashed] (0,3.5) -- (0,-.5);
\end{tpic}$$
these elements have matrices
\begin{align*}
\tau_1 = \tau_3 &= \begin{pmatrix} -t^{1/4} & 0 \\ -t^{-1/4} & t^{-3/4}
    \end{pmatrix} \\
\tau_2 &= \begin{pmatrix} t^{-3/4} & -t^{-1/4} \\ 0 & -t^{1/4}
    \end{pmatrix} \\
\sigma_3 &= \begin{pmatrix} 0 & -t^{-1} \\ 1 & -t^{-1/2} \end{pmatrix} \\
\sigma_4 &= \begin{pmatrix} 0 & t^{-3/4} \\ t^{-3/4} & 0 \end{pmatrix}.
\end{align*}
Suppose first that $t > 0$.  In this case, we use $\sigma_3$ and $\sigma_4$
as generators.  The element $\sigma_3$ is a rotation by $2\pi/3$, while the
element $\sigma_4$ is a reflection.  Then $\sigma_4$, $\sigma_3 \sigma_4
\sigma_3^{-1}$, and $\sigma_3^{-1} \sigma_4 \sigma_3$ are reflections
through three lines that make a symmetric triangle $T$; however the triangle
may have ideal or hyperideal vertices; see \fig{f:triangle}(a,b)
for examples.

The structure of a vertex of the triangle can be determined from the element
$$[\sigma_3,\sigma_4] = \tau_1^{-1} \tau_2 = 
    \begin{pmatrix} -t^{-1} & t^{-1/2} \\ -t^{-1/2} & 1-t \end{pmatrix} $$
As written, this element has determinant 1 and trace $-t+1-t^{-1}$.  It is
equivalent to negate the trace.  When $t + t^{-1} > 3$, then the eigenvalues
of $\gamma = \tau_2^{-1} \tau_1$ are real and positive, so that $\gamma$
is hyperbolic and the vertices of $T$ are hyperideal.  In the marginal
case that $t+t^{-1} = 3$, $\gamma$ is parabolic.   Finally when $t +
t^{-1} < 3$, $\gamma$ is elliptic.  Using the formula
$$t - 1 + t^{-1} = 2\cos \theta,$$
as $t$ goes to 1 from either side, $\theta$ goes monotonically to $\pi/3$
from below.

It is easy to see that for all of the values of $t > 0$ listed as discrete,
the triangle $T$ tiles the hyperbolic plane $\H^2$ by reflections through
its sides.  $B_4$ acts on this tiling and the action is discrete.  If we take
the barycentric subdivision of each copy of $T$, the result is a tiling by
a triangle $T'$ with angles of $\pi/2$ and $\pi/3$, and either an ideal or
a hyperideal vertex or an angle of $\theta/2$.  (See \fig{f:triangle}(c).)
The fundamental domain of the action of $B_4$ is either one or two copies
of $T'$.

The more subtle fact is that the action is not discrete when $\theta \ne
2\pi/n$.  This is known from the classification of 2-dimensional orbifolds.

Now let $t < 0$.  In this case we pass to the basis
$$\begin{tpic}[style=basic,scale=1.5]
\draw (0,0) .. controls (4,0) and (4,3) .. (0,3);
\draw (0,1) .. controls (2,1) and (2,2) .. (0,2);
\draw[darkgray,dashed] (0,3.5) -- (0,-.5);
\end{tpic} \qquad
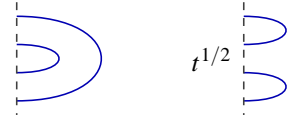
\begin{tpic}[style=basic,scale=1.5]
\draw (-.1,1.5) node[anchor=east] {$t^{1/2}$};
\draw (0,0) .. controls (2,0) and (2,1) .. (0,1);
\draw (0,2) .. controls (2,2) and (2,3) .. (0,3);
\draw[darkgray,dashed] (0,3.5) -- (0,-.5);
\end{tpic}$$
following the definition of $X(4 \cdot 1,0,t)_\R$.  Now the matrices are
\begin{align*}
\tau_1 = \tau_3 &= \begin{pmatrix} -t^{1/4} & 0 \\ -t^{-3/4} & t^{-3/4}
    \end{pmatrix} \propto
    \begin{pmatrix} -t & 0 \\ 1 & 1 \end{pmatrix} \\
\tau_2 &= \begin{pmatrix} t^{-3/4} & -t^{1/4} \\ 0 & -t^{1/4}
    \end{pmatrix} \propto
    \begin{pmatrix} 1 & -t \\ 0 & t \end{pmatrix} \\
\sigma_3 &= \begin{pmatrix} 0 & -t^{-1/2} \\ t^{-1/2} & -t^{-1/2}
    \end{pmatrix} \propto
    \begin{pmatrix} 0 & -1 \\ 1 & -1 \end{pmatrix}\\
\sigma_4 &= \begin{pmatrix} 0 & t^{-1/4} \\ t^{-5/4} & 0 \end{pmatrix}
\propto \begin{pmatrix} 0 & t \\ 1 & 0 \end{pmatrix}.
\end{align*}
After rescaling to make the matrix real, $\det \sigma_4 > 0$
and $\sigma_4$ is a rotation by $\pi$ rather than a reflection.  Meanwhile
$\sigma_3$ is a rotation by $2\pi/3$ as before.  The product $\tau_1 =
\sigma_4 \sigma_3^{-1}$ has eigenvalue ratio $-t$; it is parabolic when $t
= -1$ and is a translation otherwise.  Thus the action of $T$ preserves
the same tiling by $T'$ as before, where here $T'$ always has an ideal
or hyperideal vertex.  This time $B_4$ acts by the orientation-preserving
subgroup of the symmetry group of the tiling.

Now consider case 3.  In this case $B_4$, and in particular $\sigma_4$,
preserves the invariant, indefinite Hermitian form on $X(4 \cdot 1,0,t)$.
Thus the action of $B_4$ is contained in a copy of $\PSU(1,1) \sim
\PSL(2,\R)$, and not just in its normalizer isomorphic to $\PGL(2,\R)$.
The analysis from the previous paragraph continues, except that the
eigenvalue ratio $-t$ of $\tau_1$ now tells us that $\tau_1$ is a rotation
by $\pi-|\theta|$.  Again, when this is of the form $2\pi/n$, the action
of $B_4$ preserves a tiling of $T'$, where $T'$ is the finite triangle with
angles of $\pi/2$, $\pi/3$, and $(\pi-|\theta|)/2$.  By the classification
of 2-dimensional orbifolds, the action of $B_4$ is indiscrete for other
choices of $\theta$.

Finally in case 4, we once again use the element
$$g = [\tau_2,\tau_2\tau_3^3\tau_2\tau_1^{-1}]$$
with characteristic polynomial
$$\chi_g(x) = (x-1)(x^2+(t-1+t^{-1})^3 x + 1).$$
When
$$t + t^{-1} = 1 + 2\cos \frac{2\pi}{7},$$
the roots $\lambda$ other than $1$ have the form $\exp(i \theta)$ with
$$|\theta| \approx 165.812896^\circ.$$
This time, $\lambda$ has degree at most 6 over $\Q$.  If it did have
order $s$, then $\phi(s) \le 6$ so that $s \le 18$.  The approximation to
$\theta$ shows that $\lambda$ is not a root of unity and $g$ is elliptic
with infinite order.  Thus by \lem{l:rot}, the action of $B_4$ is indiscrete.
\end{proof}

\begin{remark} When $t = -1$ or $t = \exp(\pm 2\pi i/3)$, the action
of $B_4$ on the space $W(4 \cdot 1,0,t)$ gives a natural extension of
\thm{th:discrete} because it is always 2-dimensional.  When $t
= -1$, the triangle $T'$ has an ideal vertex and the projective action of
$B_4$ is equivalent to $\PSL(2,\Z)$.  When $t = \exp(\pm 2\pi i/3)$, then
the action is equivalent to the orientation-preserving symmetries of the
standard tiling of the Euclidean plane by equilateral triangles. Finally
when $t = 1$, the action is finite and equivalent to the symmetry group
of a single triangle.
\end{remark}

\begin{proof}[Proof of \cor{c:tutte}] It is well-known that the weighted
Potts model of a planar graph $G$ can be realized within the Kauffman
bracket by a skein replacement of every edge.  Let $Z(G,n)$ be the Potts
model with $n$ colors.  We first let
$$n = [2]^2 = t + 2 + t^{-1}.$$
Then we can make a Kauffman skein $S$ from $G$ so that
$$Z(G,n) = (-[2])^v \braket{S},$$
where $G$ has $v$ vertices.  The skein $S$ is obtained by replacing each
vertex by doubling each incoming edge to two arcs, and stitching together
these arcs into a planar matching:
$$\begin{tpic}[style=basic]
\draw (-2,0) -- (2,0); \draw (0,-2) -- (0,2);
\fill (0,0) circle (.3);
\end{tpic} \;=\;
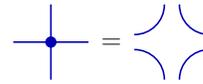
\begin{tpic}[style=basic]
\draw (-2,.4) arc (-90:0:1.6); \draw (.4,2) arc (180:270:1.6);
\draw (2,-.4) arc (90:180:1.6); \draw (-.4,-2) arc (0:90:1.6);
\end{tpic}$$
Then if an edge has Potts weight $y$, we can replace it as follows:
e can replace an
edge with Potts weight $y$ by
\eq{e:replace}{
\begin{tpic}[style=basic]
\draw (-2,0) -- (2,0);
\fill (-2,0) circle (.3); \fill (2,0) circle (.3);
\draw (0,.1) node[anchor=south] {$y$};
\end{tpic}
\;=\;
\frac{1-y}{[2]}
\begin{tpic}[style=basic]
\draw (-1,-1) arc (-45:45:1.414); \draw (1,1) arc (135:225:1.414);
\end{tpic}
\;+\;
\begin{tpic}[style=basic]
\draw (1,-1) arc (45:135:1.414); \draw (-1,1) arc (225:315:1.414);
\end{tpic}.}

In the case that $q > 4$, we can take $t > 1$ and directly apply
case 4 of \cor{c:real}, because a crossing is proportional to a weighted
edge with a real weight.  This tells us that the action of the
edge operators is real Zariski dense in $\PSL(V(n)_\R)$.  At the same time, the
action includes a Lie group of positive dimension, because the edge
operator $A_{j,y}$ has a free parameter $y$.  Thus the action 
includes all of $\PSL(V(n)_\R)$.

The case $q = 4$ is more of a corollary of the proof of \thm{th:complex}
and case 4 of \cor{c:real}.  First, even though the corresponding value of
the Jones polynomial is trivial, the Kauffman bracket still exists when
$t = 1$.  Instead of using crossings, the proofs still hold using the
replaced edge operator in equation~\eqref{e:replace}.  This time, instead
of an eigenvalue ratio of $t^{c+1}$, the eigenvalue ratio is unrestricted
because $y$ is a free parameter.  And the action includes a Lie group of
positive dimension, so the action includes all of $\PSL(V(n)_\R)$.
\end{proof}

% \bibliography{qp,qa,me,shared}

\providecommand{\bysame}{\leavevmode\hbox to3em{\hrulefill}\thinspace}
\providecommand{\MR}{\relax\ifhmode\unskip\space\fi MR }
\providecommand{\MRhref}[2]{%
  \href{http://www.ams.org/mathscinet-getitem?mr=#1}{#2}
}
\providecommand{\href}[2]{#2}

\end{document}